\documentclass[12pt]{article}
\usepackage[utf8]{inputenc}
\usepackage[english]{babel}
\usepackage{amssymb,amsmath,amsfonts,amsthm,amscd,latexsym,verbatim,graphics,epsfig,indentfirst,xcolor,comment}
\usepackage{geometry}
\usepackage{xcolor}
\def\alg{\mathrm{alg}}
\def\Var{\mathrm{Var}}

\def\RB{\mathrm{RB}}

\def\Lie{\mathrm{Lie}}
\def\Pois{\mathrm{Pois}}
\def\Com{\mathrm{Com}}
\def\NS{\mathrm{NS}}

\begin{document}

\begin{center}
{\Large Free Poisson Rota---Baxter algebra}

V. Gubarev
\end{center}

\begin{abstract}
We construct a free Poisson algebra endowed with a Rota---Baxter operator.
The same construction works for a free Poisson algebra endowed with a Nijenhuis operator.

\end{abstract}

\section{Introduction}

A linear operator $R$ defined on an algebra $A$ over a field $F$ is called a Rota---Baxter operator (RB-operator for short) of weight $\lambda \in F$ if the following relation holds:
$$
R(x)R(y) = R( R(x)y + xR(y) + \lambda xy), \quad x,y\in A.
$$
An algebra equipped with an RB-operator defined on it is called a Rota?Baxter algebra (RB-algebra for short).

The study of commutative Rota---Baxter algebras originated in a 1960 paper by G.~Baxter~\cite{Baxter60}, which addressed a specific problem in analysis. Subsequently, G.-C. Rota and others~\cite{Rota68,Cartier72} explored the combinatorial properties of RB-operators and algebras. A significant development occurred in the 1980s with the discovery of a close connection between Lie RB-algebras and solutions to the classical Yang---Baxter equation~\cite{BelaDrin82,Semenov83}. Today, RB-operators find extensive applications across mathematical physics, combina\-torics, number theory, and operad theory~\cite{Hopf03,FardThesis,GuoMonograph,GuoKeigher2000}.

In the works of G.-C. Rota, P. Cartier, and L. Guo~\cite{Rota68,Cartier72,GuoKeigher2000}, various constructions were proposed to describe free commutative RB-algebras. In 2008, K. Ebrahimi-Fard and L. Guo~\cite{FardGuo08} provided a description of free associative RB-algebras using rooted trees. Subsequently, in 2010, L.A. Bokut et al.~\cite{BokutChenDen2010} characterized these algebras through the use of Gr\"{o}bner---Shirshov bases. 
In 2016--17, different bases of free Lie RB-algebras were constructed in the works of the author~\cite{Gub2016}, the author and P.~Kolesnikov~\cite{GubKol2017}, and J.~Qiu and Y.~Chen~\cite{Chen2017}.

In 2000, M. Aguiar~\cite{Aguiar00} defined and studied Rota---Baxter operators of weight zero on Poisson algebras. By the definition, an RB-operator on a Poisson algebra has to satisfy Rota---Baxter identity with respect to both operations.
Therein, he defined a prePoisson algebra, a vector space endowed with two bilinear operations~$\circ$ and $*$ such that $\circ$ is precommutative, $*$ is preLie, and two analogues of the Leibniz identity hold.
M. Aguair showed that every RB-operator~$R$ of weight zero on a given Poisson algebra~$(A,\cdot,[,])$ defines a prePoisson algebra structure on~$A$ as follows:
$x\circ y = R(x)y$, $x*y = [R(x),y]$.
Further, different generalizations of prePoisson algebras have appeared~\cite{Aguiar20,Das,LiSun}.
In 2013, the notion of postPoisson algebra and their connection with RB-operators of nonzero weight on a Poisson algebra were established in~\cite{PoissonBialg}.
In 2013~\cite{GubKol2013}, it was proved, in particular, that every pre- and postPoisson algebra injectively embeds into its universal enveloping Poisson RB-algebra. Koszul dual structures to prePoisson algebras were studied in~\cite{Uchino}.
Pre- and postPoisson algebras appeared in~\cite{Dotsenko,DotsenkoTamaroff} as a tool to deal with operads of pre- and postLie algebras.

In this paper, we construct a basis of the free Poisson RB-algebra. 
It is known that we may construct a basis of the free Poisson algebra generated by a set~$X$ as follows.
By the Leibniz identity $[x,y\cdot z] = [x,y]\cdot z + y\cdot [x,z]$, we may assume that we initially compute all Lie products and only afterwards apply the commutative products.
Hence, one may conclude that $\Pois\langle X\rangle \cong \Com\langle \Lie\langle X\rangle\rangle$~\cite{Shestakov}.
To obtain the free Rota---Baxter Poisson algebra, we apply some similar strategy.
Now, we apply the relation 
$[a,R(b)\cdot R(c)] = [a,R(b)]\cdot R(c) + [a,R(c)]\cdot R(b)$
to avoid on the Lie algebra level the leading part of all products of the form $R(b)\cdot R(c)$.
Since we may apply operator $R$ any number of times, we have a countable numbers of steps 
of the form: RB Lie algebra level, then commutative RB level.
Roughly speaking, we construct $\RB\Pois\langle X\rangle$ inside
$\ldots \RB\Com\langle \RB\Lie\langle \ldots \RB\Com \langle\RB\Lie\langle X\rangle \rangle \ldots \rangle\rangle\ldots$
subject to certain prescribed restrictions.
We propose a potential basis and then define on the linear span of it all operations of Rota---Baxter Poisson algebra and finally prove that the resulting RB Poisson algebra is free.

We apply the obtained basis for the free Poisson algebra to construct the free NS-Poisson algebra, i.\,e. a free Poisson algebra endowed with a Nijenhuis operator. Recall that a~Nijenhuis operator defined on an algebra~$A$ satisfies the identity
$$
N(x)N(y) = N( N(x)y + xN(y) - N(xy)), \quad x,y\in A.
$$
Such operators on Lie algebras appeared in 1979~\cite{GelfandDorfman} in the context of Hamiltonian pairs.
Since 2000, Nijenhuis operators have been studied on associative algebras as well, see~\cite{Nijenhuis,FreeNijenhuisAs}.
By analogy with how one may obtain pre- and postalgebras of a variety~$\Var$ starting with 
an algebra from~$\Var$ endowed with a Rota---Baxter operator,
several authors introduced NS-associative algebras~\cite{Leroux}, NS-Lie algebras~\cite{Das} etc. defining them via Nijnehuis operators.
The general approach to defining NS-$\Var$-algebras for a given variety~$\Var$ was proposed recently~\cite{NS-algebras}.

In~\cite{BaishyaDas}, Nijenhuis operators on Poisson algebras as well as NS-Poisson algebras were defined and studied. There it was shown that the semi-classical limit of an NS-algebra deformation and a suitable filtration of an NS-algebra produce NS-Poisson algebras.

Free associative and Lie NS-algebras were constructed in~\cite{FreeNijenhuisAs} and in~\cite{Gub2016}, respectively; in both cases, they share the same basis as the free Rota---Baxter algebra of the corresponding variety. 
Thus, it is not surprising that we are able to define a free NS-Poisson algebra structure on the basis previously constructed for the free RB Poisson algebra.

The paper is organized as follows.
Section~2 provides fundamental definitions along with several examples of Poisson RB-algebras. 
Section~3 presents preliminaries on free commutative RB-algebras. 
In Section~4, we construct free Rota---Baxter Poisson algebra.
In Section~5, we apply the obtained basis to construct the free NS-Poisson algebra.

\section{Poisson Rota---Baxter algebras}

All definitions given below are well-known, see, e.\,g.~\cite{Aguiar00,BBGN}.

{\bf Definition 1}.
A vector space $V$ endowed with two bilinear products $\cdot$ and $[,]$
is called a Poisson algebra, if $(V,\cdot)$ is a commutative algebra,
$(V,[,])$ is a Lie algebra, and the following relation holds:
\begin{equation} \label{Leibniz}
[a,bc] = [a,b]c + b[a,c]
\end{equation}
for all $a,b,c\in V$.

The identity~\eqref{Leibniz} is called a Leibniz identity,
and the bracket $[,]$ is called the Poisson bracket.

{\bf Definition 2}.
Let $A$ be a Poisson algebra.
A linear operator~$R$ is called a Rota---Baxter operator (RB-operator, for short)
of weight~$\lambda$ on~$A$, if $R$ satisfies the following identities:
\begin{gather}
R(x)R(y) = R( R(x)y + xR(y) + \lambda xy), \label{eq:RBCom} \\ 
[R(x),R(y)] = R( [R(x),y] + [x,R(y)] + \lambda [x,y]). \label{eq:RBLie}
\end{gather}
A Poisson algebra with a Rota---Baxter operator is called a Poisson Rota---Baxter algebra.

{\bf Example 1}.
One may define a~Poisson bracket on the generators of the polynomial algebra $F[x,y]$
as follows: $[x,y] = y$,
extending it to the whole algebra by the Leibniz identity.
Hence, we obtain a~Poisson algebra.

We have a decomposition $F[x,y] = F[x] \oplus yF[x,y]$
into the direct vector space sum of two subalgebras. 
An operator $P$ defined on $F[x,y]$ by the rule
$$
P(f(x)) = 0,\quad
P(yg(x,y)) = - yg(x,y)
$$
is an RB-operator of weight~1 on $F[x,y]$.

Take an associative algebra $A$ and a tensor $r = \sum a_i\otimes b_i\in A\otimes A$.
Then $r$ is called a~solution to the associative Yang---Baxter equation~\cite{Aguiar00-2,Polishchuk,Zhelyabin} on $A$, if the equality
\begin{equation}\label{AYBE}
r_{13}r_{12}-r_{12}r_{23}+r_{23}r_{13} = 0
\end{equation}
holds, where
$$
r_{12} = \sum a_i\otimes b_i\otimes 1,\quad
r_{13} = \sum a_i\otimes 1\otimes b_i,\quad
r_{23} = \sum 1\otimes a_i\otimes b_i.
$$

{\bf Proposition 1}~\cite{Aguiar00}.
Let $(A,\cdot,[,])$ be a Poisson algebra and $r = \sum a_i\otimes b_i$ 
be a~solution to the associative Yang---Baxter equation on~$(A,\cdot)$. 
Then a linear operator $P\colon A\to A$ defined by the formula
$P(x) = \sum a_i x b_i$,
is a Rota---Baxter operator of weight~0 on~$A$.

Let us extend Proposition~1 to the case of nonzero weight~$\lambda\in F$.
Suppose that $A$ is an associative algebra and $r = \sum a_i\otimes b_i\in A\otimes A$.
The tensor $r$ is called a~solution to the associative Yang---Baxter equation of weight~$\lambda$~\cite{FardThesis,Ogievetsky} on~$A$, if
\begin{equation}\label{wAYBE}
r_{13}r_{12}-r_{12}r_{23}+r_{23}r_{13} = \lambda r_{13}.
\end{equation}

{\bf Proposition 2}.
Let $(A,\cdot,[,])$ be a Poisson algebra and $r = \sum a_i\otimes b_i$ 
be a~solution to the associative Yang---Baxter equation of weight~$\lambda$ on~$(A,\cdot)$.
Then a linear operator $P\colon A\to A$ defined by the formula
$P(x) = \sum a_i x b_i$,
is a Rota---Baxter operator of weight~$(-\lambda)$ on~$A$.

{\sc Proof}.
In~\cite{FardThesis}, it was shown that the operator~$P$ defined by this formula is
a~Rota---Baxter operator of weight~$(-\lambda)$ on~$(A,\cdot)$.
The proof of the remaining part that $P$ is an RB-operator of weight~$(-\lambda)$ on the algebra $(A,[,])$
follows from the commutativity of $(A,\cdot)$, the anticommutativity of~$(A,[,])$, the Leibniz identity, and the equality
$$
\sum\limits_{i,j}a_ia_jb_ib_j
 = \lambda \sum\limits_k a_kb_k
$$
holding by~\eqref{wAYBE}. 
\hfill $\square$

An algebra $A$ is called a preLie (left-symmetric), if $A$ satisfies the identity
$(x*y)*z - x*(y*z) = (y*x)*z - y*(x*z)$.
An algebra $B$ is called a precommutative (Zinbiel), if the identity 
$x(yz) = (xy+yx)z$ holds in~$A$.

{\bf Definition 3}.
A vector space $A$ endowed with two bilinear products $\circ$ and $*$
is called a~prePoisson algebra, if 
$(A,\circ)$ is a precommutative algebra,
$(A,*)$ is a preLie algebra, and the following identities hold:
\begin{gather*}
(x*y-y*x)\circ z = x*(y\circ z) - y\circ(x*z), \\
(x\circ y-y\circ x)*z = x\circ(y*z) + y\circ(x*z).
\end{gather*}

Let $A$ be a prePoisson algebra, then the vector space~$A$ considered with respect to the new operations
$x\cdot y = x\circ y + y\circ x$, $[x,y] = x*y - y*x$
is a Poisson algebra.

An algebra~$A$ endowed with four bilinear products is a~postPoisson algebra, if 
it satisifes a specific set of 13 identities, see~\cite{PoissonBialg}.
An arbitrary postPoisson algebra is also provided with the structure of a Poisson algebra.

{\bf Proposition 3}.
Let $A$ be a Poisson algebra endowed with a Rota---Baxter operator~$R$ of weight~$\lambda$.

a) \cite{Aguiar00} If $\lambda = 0$, then the operations 
$x\circ y = R(x)y$ and $x*y = [R(x),y]$ 
endow~$A$ with a~prePoisson algebra structure.

b) \cite{BBGN} If $\lambda \neq 0$, then the operations 
$x\circ y = R(x)y$, $x\cdot y = \lambda xy$, $x*y = [R(x),y]$, $x\star y = \lambda [x,y]$ 
endow~$A$ with a~postPoisson algebra structure.

In~\cite{GubKol2013}, it was proved that the embedding of any pre- and postPoisson algebra into the universal enveloping Rota---Baxter algebra 
of the corresponding weight is injective.

\section{Free commutative Rota---Baxter algebra}

Let $X$ be a non-empty set.
It is known that the set of words 
$$
R(R( R(\ldots R(R(u_n)u_{n-1})\ldots)u_3)u_2)u_1,
$$
where $u_j$, $j=1,\ldots,n-1$, are monomials in $X$ or equal to~1, $u_n$ is a~monomial in~$X$,
forms a basis of the free commutative Rota---Baxter algebra $\RB\Com\langle X\rangle$.
The product in $\RB\Com\langle X\rangle$ is defined by induction:
$$
R(a)u*R(b)v
 = R(R(a)*b + a*R(b) + \lambda a*b)uv.
$$
Let us call such basis of the free commutative RB-algebra a~standard one.

The number of appearances of the symbol $R$ in the notation of an element~$w$ of the standard basis of a free commutative RB-algebra is called $R$-degree of $w$ denoted by $\deg_R(w)$.

Suppose that a set $X$ is well-ordered.
We define the lexicographical order on the standard basis of $\RB\Com\langle X\rangle$ as follows:

1) $x<R(u)$ for every $R$-letter $R(u)$ and each $x\in X$,

2) $R(u)<R(v)$ if and only if $u<v$.

We call an element~$w$ of the standard basis of $\RB\Com\langle X\rangle$ {\it expressible}, if, ignoring the~coefficient, $w$~equals the leading word of the product $R(a)*R(b)$, where $R(a)$ and $R(b)$ are some elements of the standard basis of $\RB\Com\langle X\rangle$. 
It means that $\overline{R(a)*R(b)} = kw$, where $k\in\mathbb{N}_{>0}$ is 
an appropriate number.

For example, the word $w = R(R(x)x)$, where $x\in X$, is expressible, since
$$
R(x)*R(x) = 2R(R(x)x) + \lambda R(x^2).
$$
So, $\overline{R(x)*R(x)} = 2R(R(x)x) = 2w$.

We do not provide details about the construction of the free Rota---Baxter Lie algebra from~\cite{Gub2016}, since we actually need only the existence of such object. The notion of $R$-degree as well as order are well-defined for words forming the basis of the free RB Lie algebra.

\section{Construction of free Poisson Rota---Baxter algebra}

Let $X$ be a non-empty well-ordered set.
We define sets $V_n$, $U_n$, $W_n$ by induction on~$n$.

{\sc Step 1}. 
Define $S_1$ as the standard basis of the algebra $\RB\Lie\langle X\rangle$
and put $U_1 = S_1\cap \langle X,R(X)\rangle_{\alg}$.
By $\langle Y\rangle_{\alg}$ for given~$Y$ we mean a subalgebra of $\RB\Lie\langle X\rangle$ 
generated by~$Y$ without action of the Rota---Baxter operator.
Denote by $V_1$ the set of all monomials from $F[U_1]$ 
that contain at most one $R$-letter at the outer level, here we forget about the Lie level and $R$-letters used there.
Hereinafter, we assume that an order defined on the generators of all involved polynomial algebras satisfies the conditions~1) and~2).
Let $W_1$ be the~subset of $V_1$ consisting of all $R$-letters that are expressible words in $\RB\Com\langle U_1\rangle$.

Suppose that sets $U_s$, $V_s$, $W_s$ are defined for all numbers $s<k$.

{\sc Step $k$}. Define $S_k$ as the standard basis of the algebra 
$\RB\Lie\langle X\cup R(V_{k-1}\setminus W_{k-1})\rangle$ and put
$$
U_k = S_k\cap \langle X,R(V_{k-1}\setminus W_{k-1})\rangle_{\alg},\quad
\widetilde{U}_k = U_k\cup W_{k-1}.
$$
Denote by $V_k$ the set of all monomials from $F[\widetilde{U}_k]$
that contain at most one $R$-letter at the outer level.
For $W_k$, we take the subset of $V_k$ consisting of all $R$-letters that are 
expressible words in $\RB\Com\langle \widetilde{U}_k\rangle$.

Hence, we have defined sets $V_n$, $U_n$, $W_n$ for all~$n$.
We see that $V_n\subset V_{n+1}$, $U_n\subset U_{n+1}$, and $W_n\subset W_{n+1}$ for all~$n$. 
Finally, we define sets
$$
V_\infty = \cup V_i, \quad
U_\infty = \cup U_i, \quad
W_\infty = \cup W_i.
$$

{\bf Remark 1}.
Let us clarify why we exclude the sets $W_k$ while constructing the inter\-mediate free Lie Rota---Baxter algebra.
For this, we consider the equality that holds in any Poisson RB-algebra:
\begin{equation} \label{RBvsPoisson}
[a,R(b)*R(c)] = [a,R(b)]*R(c) + [a,R(c)]*R(b).
\end{equation}
By rewriting $R(b)*R(c)$ in the left-hand side of~\eqref{RBvsPoisson} by~\eqref{eq:RBCom},
we express $[a,\overline{R(b)*R(c)}]$ in terms of smaller words from the left-hand and the right-hand sides of~\eqref{RBvsPoisson}.

{\bf Example 2}.
Let $x,y,z,t\in X$ and $z>t$. Then

a) $R(R(x)y)z \in V_3$, since $R(x)y\in V_2$;

b) $R^3([R^2(x),y]\cdot [z,t]) \in V_6$, since 
$[R(R(x)),y],[z,t]\in W_2$ and so $[R^2(x),y]\cdot [z,t]\in V_3$.

Let $D$ be a linear span of the set $V_\infty$.
We define on~$D$ bilinear operations $\cdot$ and $[,]$ and a linear operator~$R$.

1) Let $u,v \in V_n$.
If at least one of $u$ and $v$ does not contain $R$-letters, then 
$u\cdot v$ is defined as the product of monomials in $F[\widetilde{U}_n]$.
Otherwise we present $u = R(s_1)u'$ and $v = R(s_2)v'$ and define
\begin{equation} \label{Com-product}
u\cdot v
 := R( R(s_1)\cdot s_2 + s_1\cdot R(s_2) + \lambda s_1\cdot s_2)(u'\cdot v').
\end{equation}
The product $u'\cdot v'$ is given in $F[\widetilde{U}_n]$.
The products $R(s_1)\cdot s_2,s_1\cdot R(s_2),s_1\cdot s_2$ are defined by induction on 
$\deg_R(u) + \deg_R(v)$ (in $\RB\Com\langle \widetilde{U}_n\rangle$).
Thus, we obtain a linear combination of elements from $V_t$ for some~$t$.

2) Let $u,v \in V_n$, we present $u = a_1\ldots a_k$ and $v = b_1\ldots b_l$
for $a_i,b_j\in \widetilde{U}_n$.
Given a~monomial $z = t_1\ldots t_s$, denote $z(i) = t_1\ldots t_{i-1}t_{i+1}\ldots t_s$.
We define 
\begin{equation} \label{Lie-product}
[u,v]
 = \sum\limits_{i=1}^k\sum\limits_{j=1}^l
 [a_i,b_j]\cdot( u(i)\cdot v(j) ),
\end{equation}
it remains to clarify the definition of $[a_i,b_j]$.

a) Let $a_i,b_j\not\in W_n$; then the product $[a_i,b_j]$ 
is defined in the algebra $\RB\Lie\langle X\cup R(V_{n-1}\setminus W_{n-1})\rangle$.

b) Let $a_i,b_j\in W_n$, thus, $a_i = R(c)$ and $b_j = R(d)$ for some $c,d\in V_{n-1}$.
In this case, we define
\begin{equation} \label{[ai,bj]}
[a_i,b_j] = R([R(c),d] + [c,R(d)] + \lambda [c,d]),
\end{equation}
and the products under the action of $R$ are defined by induction on the total $R$-degree of the factors.

c) Let $a_i\in W_n$ and $b_j\not\in W_n$ (in the opposite case we define $[a_i,b_j] = - [b_j,a_i]$).
Since $a_i$ is an expressible word,
it equals the leading word of the product $R(e)\cdot R(f)$
with an appropriate coefficient $k\in\mathbb{N}_{>0}$ and some $e,f\in V_{n-1}$.
Then we define
\begin{multline} \label{aiInWnbjNotInWn}
k[a_i,b_j]
 = [\overline{R(e)\cdot R(f)},b_j] \\
 := [R(e),b_j]\cdot R(f) + [R(f),b_j]\cdot R(e) - [(R(e)\cdot R(f)-ka_i),b_j].
\end{multline}
In the first factor of the product $[(R(e)\cdot R(f)-ka_i),b_j]$,
we have words with the same collection of letters as in $a_i$ but lower than $a_i$.
The products $[R(e),b_j]$ and $[R(f),b_j]$ are defined by induction on 
the total $R$-degree of the factors.

3) Let $v\in V_n$. Then we define action $R$ on $v$ as $R(v)\in V_{n+1}$.

{\bf Theorem 1}.
The space $D$ with respect to operations $\cdot,[,],R$ forms a free Poisson Rota---Baxter algebra generated by~$X$.

{\sc Proof}.
Let $u,v\in V_n$; then the product $u\cdot v$ is commutative by definition.
Let $u,v,w\in V_n$. If at least one of $u,v,w$ does not contain $R$-letters, then the identity
$(u\cdot v)\cdot w = u\cdot(v\cdot w)$ follows from the definition of~$\cdot$.
Now, consider the case where $u = R(s_1)t_1$, $v = R(s_2)t_2$, and $w = R(s_3)t_3$.
To prove associativity for such a triple of elements, we need to check the following equality:
$$
(R(s_1)\cdot R(s_2))\cdot R(s_3)
 = R(s_1)\cdot (R(s_2)\cdot R(s_3)).
$$
This relation holds by the known computations presented in the free commutative Rota---Baxter algebra~\cite{GuoMonograph}.

Thus, $D$ is a commutative RB-algebra.

Anticommutativity of the operation $[,]$ holds by definition.

Let us prove simultaneously the Leibniz identity for a triple $u,v,w$ and the Rota---Baxter identity~\eqref{eq:RBLie} for a pair $a,b$ by induction on the total $R$-degree of the factors
$\deg_R(u) + \deg_R(v) + \deg_R(w)$ and $\deg_R(a) + \deg_R(b)$.
We assume that both identities are proved when 
$\deg_R(u) + \deg_R(v) + \deg_R(w),\deg_R(a) + \deg_R(b)<r_0$
and establish them when 
$\deg_R(u) + \deg_R(v) + \deg_R(w) = \deg_R(a) + \deg_R(b) = r_0$.

We start with the Rota---Baxter identity~\eqref{eq:RBLie}.
Let $a,b\in V_n$ be two $R$-letters, i.\,.e. $a = R(a')$ and $b = R(b')$. 

{\sc Case 1}: $a,b\in U_n$. Then
\begin{equation} \label{RBLie-Proof}
[R(a'),R(b')] = R([R(a'),b'] + [a',R(b')] + \lambda [a',b'])
\end{equation}
holds since this equality holds true in the algebra 
$\RB\Lie\langle X\cup R(V_{n-1}\setminus W_{n-1})\rangle$.
This provides the induction base for~\eqref{eq:RBLie}.

{\sc Case 2}: $a,b\in W_n$. 
Then~\eqref{RBLie-Proof} is satisfied by the definition of the operation~$[,]$.

{\sc Case 3}: $a\in U_n$, $b\in W_n$.
Then $kb = \overline{R(f)\cdot R(g)}$ for some $f,g\in V_{n-1}$ and $k\in\mathbb{N}_{>0}$.
Denote $s = kb - R(f)\cdot R(g)$, $s = R(s')$. 
This implies the equality of arguments under the action of~$R$:
$kb' = s' + R(f)\cdot g+ R(g)\cdot f+ \lambda f\cdot g$.
Let us write down the expressions:
$$
kb = s + R(f)\cdot R(g), \quad
kb' = s' + R(f)\cdot g+ R(g)\cdot f+ \lambda f\cdot g.
$$
We apply the Leibniz identity for an arbitrary $d\in D$ such that $\deg_R(d) + \deg_R(b')<r_0$:
\begin{multline*}
[d,kb']
 = [d,s' + R(f)\cdot g + R(g)\cdot f + \lambda f\cdot g] \\
 = [d,s'] + [d,R(f)]\cdot g + [d,g]\cdot R(f)
 + [d,R(g)]\cdot f + [d,f]\cdot R(g)
 + \lambda([d,f]\cdot g + [d,g]\cdot f).
\end{multline*}
On the one hand, we get with the help of the Leibniz identity that
\begin{multline*}
kR([R(a'),b'] + [a',b] + \lambda [a',b']) \\
 = R( [R(a'),s'] + [R(a'),R(f)]\cdot g + [R(a'),g]\cdot R(f)
 + [R(a'),R(g)]\cdot f + [R(a'),f]\cdot R(g) \\
 + \lambda [R(a'),f]\cdot g + \lambda [R(a'),g]\cdot f
 + [a',s] + [a',R(f)]\cdot R(g) + [a',R(g)]\cdot R(f) ] \\
 + \lambda [a',s'] + \lambda [a',R(f)]\cdot g
 + \lambda [a',R(g)]\cdot f
 + \lambda [a',f]\cdot R(g)
 + \lambda [a',g]\cdot R(f)
 + \lambda^2([a',f]\cdot g + [a',g]\cdot f)).
\end{multline*}
On the other hand, applying the definition of $[,]$, we get the sum
\begin{multline*}
k[R(a'),R(b')]
 = [R(a'),s] + [R(a),R(f)]\cdot R(g) + [R(a),R(g)]\cdot R(f) \\
 = [R(a'),R(s')]
 + R([R(a'),f]+[a',R(f)]+\lambda[a',f])\cdot g \\
 + ([R(a'),f]+[a',R(f)]+\lambda[a',f])\cdot R(g)
 + \lambda([R(a'),f]+[a',R(f)]+\lambda[a',f])\cdot g \\
 + R([R(a'),g]+[a',R(g)]+\lambda[a',g])\cdot f
 + ([R(a'),g]+[a',R(g)]+\lambda[a',g])\cdot R(f) \\
 + \lambda([R(a'),g]+[a',R(g)]+\lambda[a',g])\cdot f).
\end{multline*}
To finish Case~3 it remains to check~\eqref{eq:RBLie} for the product $[R(a'),R(s')]$.
Thus, we have reduced the case $a\in U_n$, $b\in W_n$, to the one with the same $a$ and $s$ which is less than~$b$. Moreover, all expressions appearing in~$s$ have the same collection of letters as in $b$ and $R$-degree less or equal to~$r_0$. 
For all summands from $s$ with $R$-degree less than $r_0-\deg_R(a)$, we are done by the general induction.
For summands from $s$ which lie in $U_n$, we are done by Case~1.
Hence, we move on to summands from~$s$ which belong to $W_n$.
Since there is only a~finite number of words with a~fixed collection of letters and fixed $R$-degree,
the induction on the order $<$ for elements of $W_n$ is well-defined.

Now, we continue with the Leibniz identity.
Let $u = a_1\ldots a_k$, $v = b_1\ldots b_l$, $w = c_1\ldots c_m$, where $u,v,w\in V_n$.
Suppose that we have at most one $R$-letter among the letters $b_i,c_j$, then,
by the definitions of both operations $\cdot$ and $[,]$, we rewrite
$$
[u,v\cdot w]
 = \sum\limits_{i=1}^k\sum\limits_{j=1}^l[a_i,b_j]u(i)v(j)w
 + \sum\limits_{i=1}^k\sum\limits_{q=1}^m[a_i,c_q]u(i)w(q)v
 = [u,v]w + [u,w]v.
$$
Hence, we provide here the base case for the Leibniz identity.

Let $b_1 = R(b')$ and $c_1 = R(c')$. Then
\begin{multline*}
[u,v\cdot w]
 = [u,R(R(b')\cdot c' + b'\cdot R(c') + \lambda b'\cdot c')v(1)w(1)] \\
 = \sum\limits_{i=1}^k [a_i,R(R(b')\cdot c' + b'\cdot R(c') + \lambda b'\cdot c')]u(i)v(1)w(1) \\
 + \sum\limits_{i=1}^k\sum\limits_{j=2}^l[a_i,b_j]R(R(b')\cdot c' + b'\cdot R(c') + \lambda b'\cdot c')u(i)v(1,j)w(1) \\
 + \sum\limits_{i=1}^k\sum\limits_{q=2}^m[a_i,c_q]R(R(b')\cdot c' + b'\cdot R(c') + \lambda b'\cdot c')u(i)v(1)w(1,q).
\end{multline*}
Here $z(i,j)$ denotes $z_1z_2\ldots z_{i-1}z_{i+1}\ldots z_{j-1}z_{j+1}\ldots z_r$.

On the other hand,
\begin{multline*}
[u,v]w + [u,w]v
 = \sum\limits_{i=1}^k[a_i,R(b')]u(i)v(1)w
 + \sum\limits_{i=1}^k\sum\limits_{j=2}^l[a_i,b_j]u(i)v(1,j)w(1)R(b')\cdot R(c') \\
 + \sum\limits_{i=1}^k[a_i,R(c')]u(i)vw(1)
 + \sum\limits_{i=1}^k\sum\limits_{q=2}^m[a_i,c_q]u(i)v(1)w(1,q)R(b')\cdot R(c').
\end{multline*}
Note that it is necessary and sufficient to prove that 
\begin{equation} \label{(*)}
[a_i,R(R(b')\cdot c' + b'\cdot R(c') + \lambda b'\cdot c')]
 = [a_i,R(b')]\cdot R(c') + [a_i,R(c')]\cdot R(b').
\end{equation}
Hence, we may assume without loss of generality that we are proving the inductive step of the Leibniz identity for the~triple 
$a_i,b_1 = R(b'),c_1 = R(c')$.

Denote $td = \overline{R(b')\cdot R(c')} = tR(\Omega)$ for some 
$t\in\mathbb{N}_{>0}$ and $s = td - R(b')\cdot R(c') = R(s')$
such that $d>s$.
Thus, 
\begin{equation} \label{InnerUnderRBEq}
t\Omega = s' + R(b')\cdot c' + R(c')\cdot b' + \lambda b'\cdot c'.
\end{equation}

{\sc Case 1}: $a_i\not \in W_n$.
Applying the definition of the product $[,]$, we compute
\begin{multline*}
[a_i,R(b')\cdot R(c')]
 = [a_i,td - s]
 = [a_i,R(b')]\cdot R(c') + [a_i,R(c')]\cdot R(b')
 \pm [a_i,s] \\
 = [a_i,R(b')]\cdot R(c') + [a_i,R(c')]\cdot R(b'),
\end{multline*}
so~\eqref{(*)} is proved.

{\sc Case 2}: $a_i \in W_n$. Thus $a_i = R(\Psi)$ for suitable~$\Psi$, and the left-hand side of~\eqref{(*)} is rewritten by~\eqref{[ai,bj]} as follows:
\begin{multline*}
LL := [a_i,R(R(b')\cdot c' + b'\cdot R(c') + \lambda b'\cdot c')]
 = [a_i, td] - [a_i,s] \\
 = t[R(\Psi), R(\Omega)] - [a_i,s]
 = tR( [R(\Psi), \Omega] + [\Psi, R(\Omega)] + \lambda [\Psi, \Omega] ) - [R(\Psi),R(s')].
\end{multline*}
Note that
$$
tR( [R(\Psi), \Omega])
 = R( [R(\Psi), R(b')\cdot c' + R(c')\cdot b' + \lambda b'\cdot c']) + R([R(\Psi),s']),
$$

\vspace{-0.75cm}

\begin{multline*}
tR([\Psi,R(\Omega)])
 = R([\Psi,R(b')\cdot R(c')]) + R([\Psi,s]) \\
 = R\big([\Psi,R(b')]\cdot R(c') + [\Psi,R(c')]\cdot R(b')) + R([\Psi,R(s')]\big),
\end{multline*}
here we apply the Leibniz identity which holds by the induction hypothesis for words with a lower total $R$-degree.
Therefore
\begin{multline*}
LL = R( [R(\Psi), R(b')\cdot c' + R(c')\cdot b' + \lambda b'\cdot c']) \\
 + R([\Psi,R(b')]\cdot R(c') + [\Psi,R(c')]\cdot R(b')) 
 + \lambda R([\Psi,t\Omega-s']) \\
 \mathop{=}\limits^{\eqref{InnerUnderRBEq}} 
 R( [R(\Psi), R(b')\cdot c' + R(c')\cdot b' + \lambda b'\cdot c']) \\
 + R([\Psi,R(b')]\cdot R(c') + [\Psi,R(c')]\cdot R(b')) 
 + \lambda R([\Psi, R(b')\cdot c' + R(c')\cdot b' + \lambda b'\cdot c']).
\end{multline*}
Above, we have applied~\eqref{eq:RBLie} for $[R(\Psi),R(s')]$.

The right-hand side of~\eqref{(*)} equals
\begin{multline*}
RR:= [R(\Psi),R(b')]\cdot R(c') + [R(\Psi),R(c')]\cdot R(b') \\
 = R( [R(\Psi),b'] + [\Psi,R(b')] + \lambda [\Psi,b'] )\cdot R(c') \\
 + R( [R(\Psi),c'] + [\Psi,R(c')] + \lambda [\Psi,c'] )\cdot R(b') \allowdisplaybreaks \\
 = R\big( R( [R(\Psi),b'] + [\Psi,R(b')] + \lambda [\Psi,b'] )\cdot c'  \\
 + [R(\Psi),b']\cdot R(c') + [\Psi,R(b')]\cdot R(c') + \lambda [\Psi,b']\cdot R(c') \\
 + \lambda [R(\Psi),b']\cdot c' + \lambda [\Psi,R(b')]\cdot c' + \lambda^2 [\Psi,b']\cdot c' \big) \\
 + R\big( R( [R(\Psi),c'] + [\Psi,R(c')] + \lambda [\Psi,c'] )\cdot b'  \\
 + [R(\Psi),c']\cdot R(b') + [\Psi,R(c')]\cdot R(b') + \lambda [\Psi,c']\cdot R(b') \\
 + \lambda [R(\Psi),c']\cdot b' + \lambda [\Psi,R(c')]\cdot b' + \lambda^2 [\Psi,c']\cdot b' \big).
\end{multline*}
Here we apply~\eqref{eq:RBLie} to words with a lower total $R$-degree than $\deg_R(a_i) + \deg_R(b_1) + \deg_R(c_1)$. 

For the proof of the equality $LL = RR$, it remains to apply the Leibniz identity and~\eqref{eq:RBLie} inside the action of~$R$. By induction, both identities hold and their usage is correct.

Now, we prove the Jacobi identity for $[,]$, i.\,e. that the identity
$$
J(u,v,w) = [[u,v],w] + [[v,w],u] + [[w,u],v] = 0
$$
holds for all $u,v,w\in D$.

Let $u = a_1\ldots a_k$, $v = b_1\ldots b_l$, $w = c_1\ldots c_m$, where $u,v,w\in V_n$. 
Using the Leibniz identity, we rewrite
\begin{multline} \label{Jacobi-1}
[[u,v],w]
 = \sum\limits_{i=1}^k\sum\limits_{j=1}^l[[u_i,v_j]u(i)v(j),w]
 = \sum\limits_{i=1}^k\sum\limits_{j=1}^l\sum\limits_{q=1}^m[[u_i,v_j],w_q]u(i)v(j)w(q) \\
 + \sum\limits_{i,s=1,\,i\neq s}^k\sum\limits_{j=1}^l\sum\limits_{q=1}^m[u_i,v_j]\cdot [u_s,w_q]u(i,s)v(j)w(q) \allowdisplaybreaks \\
 + \sum\limits_{i=1}^k\sum\limits_{j,t=1,\,j\neq t}^l\sum\limits_{q=1}^m[u_i,v_j]\cdot [v_t,w_q]u(i)v(j,t)w(q).
\end{multline}

Dealing analogously with $[[v,w],u]$ and $[[w,u],v]$, 
we obtain after collecting like terms that
$$
J(u,v,w)
 = \sum\limits_{i=1}^k\sum\limits_{j=1}^l\sum\limits_{q=1}^m J(u_i,v_j,w_q)u(i)v(j)w(q).
$$
Therefore it is necessary and sufficient to prove the Jacobi identity for the letters $a,b,c\in V_n$.

{\sc Case 1}: $a,b,c\in U_n$. 
Then $J(a,b,c) = 0$, since the Jacobi identity holds in the Lie algebra 
$\RB\Lie\langle X\cup R(V_{n-1}\setminus W_{n-1})\rangle$.

{\sc Case 2}: $a,b\in U_n$ and $c\in W_n$.
Since $c$ is an expressible word, we have 
$kc = \overline{R(d)\cdot R(e)}$ for suitable $d,e\in V_{n-1}$ and $k\in\mathbb{N}_{>0}$.
Denote by $s = kc - R(d)\cdot R(e)$ the sum of $R$-letters of the same collection of letters as $c$ but smaller than $c$.
Applying the definition of the product $[,]$ for $[[b,c],a]$ and $[[c,a],b]$
and the Leibniz identity for $[[a,b],c]$ we obtain
\begin{multline*}
k([[a,b],c] + [[b,c],a] + [[c,a],b]) \allowdisplaybreaks \\
 = [[a,b],R(d)]\cdot R(e) + [[a,b],R(e)]\cdot R(d)
 + [[b,R(d)]\cdot R(e),a] + [[b,R(e)]\cdot R(d),a] \\
 + [[R(d),a]\cdot R(e),b] + [[R(e),a]\cdot R(d),b] + J(a,b,s)  \\
 = [[a,b],R(d)]\cdot R(e) + [[a,b],R(e)]\cdot R(d)
 + [[b,R(d)],a]\cdot R(e) + [b,R(d)]\cdot [R(e),a] \\
 + [[b,R(e)],a]\cdot R(d) + [b,R(e)]\cdot [R(d),a]
 + [[R(d),a],b]\cdot R(e) + [R(d),a]\cdot [R(e),b] \allowdisplaybreaks \\
 + [[R(e),a],b]\cdot R(d) + [R(e),a]\cdot [R(d),b] + J(a,b,s) \\
 = J(a,b,R(d))\cdot R(e) + J(a,b,R(e))\cdot R(d) + J(a,b,s).
\end{multline*}

The equalities 
$J(a,b,R(d)) = J(a,b,R(e)) = 0$ follow by induction on the total $R$-degree of the triple of elements and $J(a,b,s) = 0$ holds by induction on the order in~$D$. Note that $s$ equals the sum of words with the same collection of letters as in~$c$ but smaller with respect to the order~$<$ than~$c$.
Since there is only a~finite number of words with a~fixed collection of letters,
the induction on the order $<$ is correct.

{\sc Case 3}: $a\in U_n$, $b,c\in W_n$.
Then $kb = \overline{R(d)\cdot R(e)}$ and 
$lc = \overline{R(f)\cdot R(g)}$ for some $d,e,f,g\in V_{n-1}$ and $k,l\in\mathbb{N}_{>0}$.
Denote $s = kb - R(d)\cdot R(e)$ and 
$t = lc - R(f)\cdot R(g)$.
Applying the Leibniz identity, we write down
\begin{multline*}
kl[[a,b],c]
 = l([[a,s],c]
 + [[a,R(d)]\cdot R(e),c]
 + [[a,R(e)]\cdot R(d),c]) \\
 = l([[a,s],c]
 + [[a,R(d)],c]\cdot R(e)
 + [[a,R(e)],c]\cdot R(d)
 + [a,R(d)]\cdot[R(e),c]
 + [a,R(e)]\cdot[R(d),c] ) \\
 = [[a,s],t] + [[a,s],R(f)]\cdot R(g) + [[a,s],R(g)]\cdot R(f) \\
 + [[a,R(d)],t]\cdot R(e)
 + [[a,R(d)],R(f)]\cdot R(e)R(g)
 + [[a,R(d)],R(g)]\cdot R(e)R(f) \\
 + [[a,R(e)],t]\cdot R(d)
 + [[a,R(e)],R(f)]\cdot R(d)R(g)
 + [[a,R(e)],R(g)]\cdot R(d)R(f) \\
 + [a,R(d)]\cdot[R(e),t]
 + [a,R(d)]\cdot[R(e),R(f)]\cdot R(g)
 + [a,R(d)]\cdot[R(e),R(g)]\cdot R(f) \\
 + [a,R(e)]\cdot[R(d),t]
 + [a,R(e)]\cdot[R(d),R(f)]\cdot R(g)
 + [a,R(e)]\cdot[R(d),R(g)]\cdot R(f).
\end{multline*}
Analogously, we get
\begin{multline*}
kl[[a,c],b]
 = [[a,t],s] + [[a,t],R(d)]\cdot R(e) + [[a,t],R(e)]\cdot R(d) \\
 + [[a,R(f)],s]\cdot R(g)
 + [[a,R(f)],R(d)]\cdot R(e)R(g)
 + [[a,R(f)],R(e)]\cdot R(d)R(g) \\
 + [[a,R(g)],s]\cdot R(f)
 + [[a,R(g)],R(d)]\cdot R(e)R(f)
 + [[a,R(g)],R(e)]\cdot R(d)R(f) \\ \allowdisplaybreaks
 + [a,R(f)]\cdot[R(g),s]
 + [a,R(f)]\cdot[R(g),R(d)]\cdot R(e)
 + [a,R(f)]\cdot[R(g),R(e)]\cdot R(d) \\
 + [a,R(g)]\cdot[R(f),s]
 + [a,R(g)]\cdot[R(f),R(d)]\cdot R(e)
 + [a,R(g)]\cdot[R(f),R(e)]\cdot R(d).
\end{multline*}
Finally, with the help of the Lebniz identity we rewrite~$[[b,c],a]$:
\begin{multline*}
kl[[b,c],a]
 = k( [[b,t],a]
 + [[b,R(f)]\cdot R(g),a]
 + [[b,R(g)]\cdot R(f),a] ) \\
 {=} k( [[b,t],a]
 + [[b,R(f)],a]\cdot R(g)
 + [b,R(f)]\cdot [R(g),a]
 + [[b,R(g)],a]\cdot R(f)
 + [b,R(g)]\cdot [R(f),a] ) \\
 = [[s,t],a]
 + [[R(d),t],a]\cdot R(e)
 + [R(d),t]\cdot [R(e),a]
 + [[R(e),t],a]\cdot R(d)
 + [R(e),t]\cdot [R(d),a] \\
 + [[s,R(f)],a]\cdot R(g)
 + [[R(d),R(f)],a]\cdot R(e)R(g)
 + [R(d),R(f)]\cdot [R(e),a]\cdot R(g) \\
 + [[R(e),R(f)],a]\cdot R(d)R(g)
 + [R(e),R(f)]\cdot [R(d),a]\cdot R(g) \\
 + [s,R(f)]\cdot [R(g),a]
 + [R(d),R(f)]\cdot [R(g),a]\cdot R(e)
 + [R(e),R(f)]\cdot [R(g),a]\cdot R(d) \\
 + [[s,R(g)],a]\cdot R(f)
 + [[R(d),R(g)],a]\cdot R(e)R(f)
 + [R(d),R(g)]\cdot [R(e),a]\cdot R(f) \\
 + [[R(e),R(g)],a]\cdot R(d)R(f)
 + [R(e),R(g)]\cdot [R(d),a]\cdot R(f) \\
 + [s,R(g)]\cdot [R(f),a]
 + [R(d),R(g)]\cdot [R(f),a]\cdot R(e)
 + [R(e),R(g)]\cdot [R(f),a]\cdot R(d).
\end{multline*}

To compute $kl( [[a,b],c] - [[a,c],b] + [[b,c],a] )$
we use the induction and the Jacobi identity for smaller triples of basic elements. Hence, we prove that $J(a,b,c) = 0$.

{\sc Case 4}: $a,b,c\in W_n$.
Denote $a = R(a')$, $b = R(b')$, $c = R(c')$. We compute by the definition of~$[,]$ and~\eqref{eq:RBLie}
\begin{multline*}
[[R(a'),R(b')],R(c')]
 = [ R( [R(a'),b'] + [a',R(b')] + \lambda [a',b'] ) ,R(c')] \\
 = R\big( [ [R(a'),R(b')], c'] 
 + [[R(a'),b'],R(c')] + [[a',R(b')],R(c')] + \lambda [[a',b'],R(c')]  \\
 + \lambda[[R(a'),b'],c'] + \lambda[[a',R(b')],c'] + \lambda^2[[a',b'],c'] \big).
\end{multline*}
Dealing analogously wth 
$[[R(b'),R(c')],R(a')]$ and $[[R(c'),R(a')],R(b')]$,
we reduce the check of the Jacobi identity to words with a lower total $R$-degree.

Therefore we have proved that $D$ under the mentioned operations is a Poisson Rota---Baxter algebra.
Our goal is to show that $D$ is a free Poisson Rota---Baxter algebra generated by~$X$.

{\bf Theorem 2}.
We have $D\cong \RB\mathrm{Pois}\langle X\rangle$.

{\sc Proof}.
Taking into account realtion~\eqref{RBvsPoisson}, every element of $\RB\mathrm{Pois}\langle X\rangle$
can be linearly expressed in terms of words from $V_\infty$.
By the construction and Theorem~1, $D$ is a~Rota---Baxter algebra generated by~$X$. 
Hence, $D$ is a homomorphic image of the algebra $\RB\mathrm{Pois}\langle X\rangle$,
let us denote by $\varphi$ a corresponding epimorphism acting from $\RB\mathrm{Pois}\langle X\rangle$ onto~$D$. 
Note that $\ker\varphi = (0)$, otherwise~$\varphi$ would not be surjective. 
Therefore $D$ is the free Poisson RB-algebra.
\hfill $\square$

{\bf Problem}.
Can we describe the expressions obtained in the free Poisson Rota---Baxter algebra generated by a~set~$X$ starting from~$X$ and repeatedly applying the following operations:
$$
a\circ b = R(a)b, \quad
a\star b = [R(a),b]?
$$
A linear span of such expressions forms a free prePoisson algebra generated by~$X$ when the weight of the RB-operator equals zero or a free postPoisson algebra generated by~$X$ otherwise.

\section{Construction of free Poisson Nijenhuis algebra}

{\bf Definition 4}.
Let $A$ be a Poisson algebra.
A linear operator~$N$ is called a Nijenhuis operator on~$A$ if $N$ satisfies the following identities:
\begin{gather}
N(x)N(y) = N( N(x)y + xN(y) - N(xy) ), \label{eq:NSCom} \\ 
[N(x),N(y)] = N( [N(x),y] + [x,N(y)] - N([x,y]) ). \label{eq:NSLie}
\end{gather}
A Poisson algebra equipped with a Nijenhuis operator is called a Poisson Nijenhuis algebra.

Let us explain how to obtain an extension of prePoisson algebras using Nijenhuis operators.
For this, we need several definitions.

A vector space $A$ with two bilinear operations $\cdot,\wedge$ is called NS-Lie, if $\wedge$ is anticommutative
and $A$ satisfies the identities
$$
(a\star b)\cdot c = a\cdot (b\cdot c) - b\cdot (a\cdot c), \quad
(a\star b)\star c + (b\star c)\star a + (c\star a)\star b = 0,
$$
where $a\star b = a\cdot b - b\cdot a + a\wedge b$.
A vector space $B$ with two bilinear operations $\circ,\vee$ is called NS-commutative if $\vee$ is commutative and $B$ satisfies the identities
$$
a\circ (b\circ c) = (a*b)\circ c, \quad
a* (b* c) = (a*b)* c,
$$
where $a*b = a\circ b + b\circ a + a\vee b$.

{\bf Definition 5}.
A vector space $A$ endowed with four bilinear products $\cdot,\wedge,\circ,\vee$
is called a~NS-Poisson algebra if 
$(A,\circ,\vee)$ is an NS-commutative-algebra,
$(A,\cdot,\wedge)$ is an NS-Lie algebra, and the following identities hold:
\begin{gather*}
(a*b)\star c = a*(b \star c) + b*(a \star c), \quad
(a*b)\cdot c = a\circ (b\cdot c) + b\circ(a\cdot c), \\
(a\star b)\circ c = a\circ (b\cdot c) - b\cdot (a\circ c),
\end{gather*}
where $a\star b = a\cdot b - b\cdot a + a\wedge b$ and 
$a*b = a\circ b + b\circ a + a\vee b$.

{\bf Remark 2}.
The main procedure how one may write down identities of $\NS$-$\Var$-algebras for a given variety $\Var$ of algebras was presented in~\cite{NS-algebras}. 
Actually, this procedure represents the property that
$\NS$-$\Var = \NS$-$\Lie\bullet \Var$, where 
$\bullet$ is the Manin black product of the corresponding operads.

If $(A,\cdot,\wedge,\circ,\vee)$ is an NS-Poisson algebra, then the vector space~$(A,*,\star)$ is a Poisson algebra. If $\wedge$ and $\vee$ are trivial operations, i.\,e. identically equal to zero, then
such an NS-Poisson algebra is a prePoisson algebra.

{\bf Proposition 4}.
Let $(A,\perp,[,])$ be a Poisson algebra endowed with a Nijenhuis operator~$N$.
Define new products on $A$ as follows:
$$
a\circ b = N(a)\perp b, \quad 
a\vee b = - N(a\perp b), \quad
a\cdot b = [N(a),b], \quad 
a\wedge b = - N([a,b]).
$$
Then $(A,\cdot,\wedge,\circ,\vee)$ is an NS-Poisson algebra.

Let $X$ be a non-empty well-ordered set.
We can repeat the construction of the sets $V_n,U_n,W_n$ and $V_\infty,U_\infty,W_\infty$
by replacing free Rota---Baxter Lie and commutative algebras with
free Nijenhuis Lie and commutative algebras, respectively. 
Such a~replacement does not change the resulting set $V_\infty$, since 
the bases of free RB Lie and commutative algebras can be taken as the bases for the corresponding free Nijenhuis algebras~\cite{Gub2016}; these algebras differ only in their operations.

Let $D$ be the linear span of the set $V_\infty$.
We define bilinear operations $\cdot$ and $[,]$ and a linear operator~$N$ on~$D$ 
in much the same way as above, with one change: we replace 
$\lambda xy$ and $\lambda [x,y]$ with $-N(xy)$ and $-N([x,y])$, respectively.

{\bf Theorem 2}.
The space $D$ with respect to operations $\cdot,[,],N$ forms a free Nijenhuis Poisson algebra generated by~$X$.

{\sc Proof}.
We highlight only the parts that differ from the proof of Theorem~1.

{\bf I}.
The simultaneous proof of the Leibniz identity and the Lie Nijenhuis identity goes by induction on the total $N$-degree.
We begin with the Nijenhuis identity~\eqref{eq:NSLie}.
Let $a,b\in V_n$ be two $N$-letters, i.\,e. $a = N(a')$ and $b = N(b')$. 

{\sc Case 3}: $a\in U_n$, $b\in W_n$.
Then $kb = \overline{N(f)\cdot N(g)}$ for some $f,g\in V_{n-1}$ and $k\in\mathbb{N}_{>0}$.
Denote $s = kb - N(f)\cdot N(g)$, $s = N(s')$. 
This implies the following equality for the arguments under the action of~$N$:
$kb' = s' + N(f)\cdot g+ N(g)\cdot f - N(f\cdot g)$.
Let us write down the expressions:
$$
kb = s + N(f)\cdot N(g), \quad
kb' = s' + N(f)\cdot g+ N(g)\cdot f - N(f\cdot g).
$$
We apply the Leibniz identity for an arbitrary $d\in D$ such that $\deg_N(d) + \deg_N(b')<r_0$:
\begin{multline*}
[d,kb']
 = [d,s' + N(f)\cdot g + N(g)\cdot f - N(f\cdot g)] \\
 = [d,s'] + [d,N(f)]\cdot g + [d,g]\cdot N(f)
 + [d,N(g)]\cdot f + [d,f]\cdot N(g)
 - [d,N(f\cdot g)].
\end{multline*}
On the one hand, we get with the help of the Leibniz identity that
\begin{multline*}
RR = kN([N(a'),b'] + [a',N(b')] - N([a',b'])) \\
 = N\big( [N(a'),s'] + [N(a'),N(f)]\cdot g + [N(a'),g]\cdot N(f)
 + [N(a'),N(g)]\cdot f + [N(a'),f]\cdot N(g) \\
 - [N(a'),N(f\cdot g)]
 + [a',N(s')] + [a',N(f)]\cdot N(g) + [a',N(g)]\cdot N(f) \\
 - N([a',s' + N(f)\cdot g + N(g)\cdot f - N(f\cdot g)] \big).
\end{multline*}
On the other hand, applying the definition of $[,]$, we get the sum
\begin{multline*}
LL = k[N(a'),N(b')]
 = [N(a'),N(s')] + [N(a'),N(f)]\cdot N(g) + [N(a'),N(g)]\cdot N(f) \\
 = N\big( [N(a'),s'] + [a',N(s')] - N([a',s']) 
 + [N(a'),N(f)]\cdot g \\
 + [N(a'),f]\cdot N(g) + [a',N(f)]\cdot N(g) - N([a',f])\cdot N(g) 
 - N( [N(a'),f]\cdot g +[a',N(f)]\cdot g  \\
 - N([a',f])\cdot g ) 
 + [N(a'),N(g)]\cdot f + [N(a'),g]\cdot N(f) + [a',N(g)]\cdot N(f) - N([a',g])\cdot N(f) \\
 - N( [N(a'),g]\cdot f +[a',N(g)]\cdot f  - N([a',g])\cdot f) \big).
\end{multline*}

By collecting like terms we deduce that 
$LL - RR = N^2(\Delta)$ for 
\begin{multline*}
\Delta 
 = [a',N(f)]\cdot g + [a',g]\cdot N(f) + [a',f]\cdot N(g) + [a',N(g)]\cdot f - [a',N(f\cdot g)] \\
 + [N(a'),f\cdot g] + [a',N(f\cdot g)] - N([a',f\cdot g]) 
 - N([a',f])\cdot g + [a',f]\cdot N(g) - N([a',f]\cdot g) \\
 - [N(a'),f]\cdot g - [a',N(f)]\cdot g + N([a',f])\cdot g
 - N([a',g])\cdot f + [a',g]\cdot N(f) - N([a',g]\cdot f) \\
 - [N(a'),g]\cdot f - [a',N(g)]\cdot f + N([a',g])\cdot f) 
 = 0.
\end{multline*}

{\bf II}.
Next, we proceed with the Leibniz identity. It is necessary and sufficient to prove that:
\begin{equation} \label{(**)}
[a_i,N(N(b')\cdot c' + b'\cdot N(c') - N(b'\cdot c'))]
 = [a_i,N(b')]\cdot N(c') + [a_i,N(c')]\cdot N(b').
\end{equation}
Hence, without loss of generality, we may assume we are proving the inductive step for the triple 
$a_i,b_1 = N(b'),c_1 = N(c')$.

Denote $td = \overline{N(b')\cdot N(c')} = tN(\Omega)$ for some 
$t\in\mathbb{N}_{>0}$ and $s = td - N(b')\cdot N(c') = N(s')$
such that $d>s$.
Thus, 
\begin{equation} \label{InnerUnderNSEq}
t\Omega = s' + N(b')\cdot c' + N(c')\cdot b' - N( b'\cdot c' ).
\end{equation}

{\sc Case 2}: $a_i \in W_n$. Thus $a_i = N(\Psi)$ for suitable~$\Psi$, and the left-hand side of~\eqref{(**)} is rewritten as follows:
\begin{multline*}
LL := [a_i,N(N(b')\cdot c' + b'\cdot N(c') - N( b'\cdot c'))]
 = [a_i, td] - [a_i,s] \\
 = t[N(\Psi), N(\Omega)] - [a_i,s]
 = tN( [N(\Psi), \Omega] + [\Psi, N(\Omega)] - N( [\Psi, \Omega] ) ) - [N(\Psi),N(s')].
\end{multline*}
Note that
$$
tN( [N(\Psi), \Omega])
 = N( [N(\Psi), N(b')\cdot c' + N(c')\cdot b' - N( b'\cdot c'] ) ) + N([N(\Psi),s']),
$$

\vspace{-0.75cm}

\begin{multline*}
tN([\Psi,N(\Omega)])
 = N([\Psi,N(b')\cdot N(c')]) + N([\Psi,s]) \\
 = N\big([\Psi,N(b')]\cdot N(c') + [\Psi,N(c')]\cdot N(b')) + N([\Psi,N(s')]\big),
\end{multline*}
here we apply the Leibniz identity which holds by the induction hypothesis for words with a lower total $N$-degree.
Therefore
\begin{multline*}
LL = N( [N(\Psi), N(b')\cdot c' + N(c')\cdot b' - N( b'\cdot c')]) \\
 + N([\Psi,N(b')]\cdot N(c') + [\Psi,N(c')]\cdot N(b')) 
 -  N^2([\Psi,t\Omega-s']) \\
 \mathop{=}\limits^{\eqref{InnerUnderNSEq}} 
 N( [N(\Psi), N(b')\cdot c' + N(c')\cdot b' - N( b'\cdot c')]) \\
 + N([\Psi,N(b')]\cdot N(c') + [\Psi,N(c')]\cdot N(b')) 
 - N^2([\Psi, N(b')\cdot c' + N(c')\cdot b' - N(b'\cdot c')]).
\end{multline*}
Above, we have applied~\eqref{eq:NSLie} for $[N(\Psi),N(s')]$.

The right-hand side of~\eqref{(**)} equals
\begin{multline*}
RR:= [N(\Psi),N(b')]\cdot N(c') + [N(\Psi),N(c')]\cdot N(b') \\
 = N( [N(\Psi),b'] + [\Psi,N(b')] - N( [\Psi,b']) )\cdot N(c') 
 + N( [N(\Psi),c'] + [\Psi,N(c')] - N( [\Psi,c']) )\cdot N(b') \\
 = N\big( N( [N(\Psi),b'] + [\Psi,N(b')] - N( [\Psi,b']) )\cdot c'  
 + [N(\Psi),b']\cdot N(c') + [\Psi,N(b')]\cdot N(c') \\
 - N( [\Psi,b'])\cdot N(c') 
 - N([N(\Psi),b']\cdot c') - N([\Psi,N(b')]\cdot c') + N( N( [\Psi,b'])\cdot c') \big) \\
N\big( N( [N(\Psi),c'] + [\Psi,N(c')] - N( [\Psi,c']) )\cdot b'  
 + [N(\Psi),c']\cdot N(b') + [\Psi,N(c')]\cdot N(b') \\
 - N( [\Psi,c'])\cdot N(b') 
 - N([N(\Psi),c']\cdot b') - N([\Psi,N(c')]\cdot b') + N( N( [\Psi,c'])\cdot b') \big).
\end{multline*}
Here we apply~\eqref{eq:NSLie} to words with a lower total $N$-degree than $\deg_N(a_i) + \deg_N(b_1) + \deg_N(c_1)$. 

For the proof of the equality $LL = RR$, it remains to apply the Leibniz identity and~\eqref{eq:NSLie} inside the action of~$N$. 

{\bf III}.
Now, we prove the Jacobi identity for $[,]$.
{\sc Case 4}: $a,b,c\in W_n$.
Denote $a = N(a')$, $b = N(b')$, $c = N(c')$. We compute by the definition of~$[,]$ and~\eqref{eq:NSLie}
\begin{multline*}
[[N(a'),N(b')],N(c')]
 = [ N( [N(a'),b'] + [a',N(b')] - N( [a',b']) ) ,N(c')] \\
 = N\big( [ [N(a'),N(b')], c'] 
 + [[N(a'),b'],N(c')] + [[a',N(b')],N(c')] - [N( [a',b'] ),N(c')]  \\
 - N( [[N(a'),b'],c'] + [[a',N(b')],c'] - N( [N([a',b']),c']) ) \big) \\
 = N\big( [ [N(a'),N(b')], c'] + [[N(a'),b'],N(c')] + [[a',N(b')],N(c')]  \\
 - N( [[N(a'),b'],c'] + [[a',N(b')],c'] + [[a',b'],N(c')] - N([[a',b'],c'] ) \big).
\end{multline*}
Dealing analogously wth 
$[[N(b'),N(c')],N(a')]$ and $[[N(c'),N(a')],N(b')]$,
we reduce the check of the Jacobi identity to words with a lower total $N$-degree.

\section*{Acknowledgements}

The project is supported by the Russian Science Foundation (project 25-41-00005).

\noindent Vsevolod Gubarev \\
Sobolev Institute of Mathematics \\
Acad. Koptyug ave. 4, 630090 Novosibirsk, Russia \\
Novosibirsk State University \\
Pirogova str. 2, 630090 Novosibirsk, Russia \\
e-mail: wsewolod89@gmail.com

\end{document}